\def\zfc{{\rm ZFC}_0}

\newcommand{\qst}{Q^{St}}
\newcommand{\qsts}[1]{Q^{{St}}_{#1}}

\def\fol{L}

\def\st{L^{St}}
\def\sti{L_{\omega_1\omega}^{St}}

\def\ss{\langle S_0,S_1,...\rangle}
\def\ssn{\langle S_0,S_1,...,S_n\rangle}
\def\ttn{\langle T_0,T_1,...,T_n\rangle}
\def\sss{S_0,S_1,...}

\def\sst{\langle S^*_0,S^*_1,...\rangle}
\def\ssst{S^*_0,S^*_1,...}
\def\Ax{\mbox{Ax}}
\def\Val{\mbox{Val}}
\def\sd{S_{\diamondsuit}}

\newtheorem{thm}{Theorem}
\newtheorem{lem}[thm]{Lemma}
\newtheorem{cor}[thm]{Corollary}
\newtheorem{dfn}[thm]{Definition}
\newtheorem{pro}[thm]{Proposition}

\newtheorem{theorem}{Theorem}
\newtheorem{lemma}[theorem]{Lemma}
\newtheorem{corollary}[theorem]{Corollary}
\newtheorem{definition}[theorem]{Definition}
\newtheorem{proposition}[theorem]{Proposition}
\newtheorem{example}[theorem]{Example}

\newcommand{\mn}{{\mathfrak N}}
\newcommand{\mm}{{\mathfrak M}}
\newcommand{\ma}{{\mathfrak A}}

\def\F{{\cal F}}
\renewcommand{\L}{{\cal L}}
\newcommand{\I}{{\cal I}}
\newcommand{\J}{{\cal J}}
\newcommand{\Ifin}{{\cal I}_{<\oo}}

\def\proof{\noindent {\bf Proof.}\hspace{1mm} }
\def\qed{$\Box$\medbreak }

\renewcommand{\aa}{\alpha}
\newcommand{\bb}{\beta}

\newcommand{\oo}{\omega}

\renewcommand{\ll}{\lambda}
\newcommand{\kk}{\kappa}

\newcommand{\tt}{\tau}

\newcommand{\pl}{\lambda}
\newcommand{\pk}{\kappa}

\newcommand{\ps}{\sigma}

\newcommand{\al}{\aleph}
\newcommand{\ri}{\rightarrow}
\newcommand{\qmm}{Q^{\mbox{\tiny MM}}}

\newcommand{\se}{\subseteq}

\newcommand{\A}{{\cal A}}
\renewcommand{\L}{{\cal L}}

\newcommand{\open}{\mathbb}

\newcommand{\oF}{{\open F}}

\newcommand{\II}{\exists}
\newcommand{\la}{\langle}
\newcommand{\ra}{\rangle}

\renewcommand{\L}{{\cal L}}

\def\II{\la (D_n)_{n<m},(I_{n,a})_{n<m,a\in D_n},(h_n)_{n<m}\ra}
\def\IIP{\la (D'_n)_{n<m'},(I'_{n,a})_{n<m',a\in D'_n},(h'_n)_{n<m'}\ra}

\def\lal{\aleph_{\aa_l}}
\def\lan{\aleph_{\aa_n}}
\newcommand{\vecz}[1]{\mbox{\boldmath{$#1$}}}

\documentclass[12pt]{article}
\usepackage{amssymb,latexsym}
\author{
Saharon Shelah
\thanks{Research partially supported by the United States-Israel
        Binational Science
Foundation. Publication number [ShVa:790] }\\
Institute of Mathematics\\
Hebrew University\\
Jerusalem, Israel\\
\and
Jouko V\"a\"an\"anen
\thanks{Research partially supported by
grant 40734 of the Academy of Finland.}\\
Department of Mathematics\\
University of Helsinki\\
Helsinki, Finland}
\title{Recursive Logic Frames}
\begin{document}

\maketitle

\begin{abstract}
We define the concept of a {\em logic frame}, which extends the
concept of an abstract logic by adding the concept of a syntax and
an axiom system. In a {\em recursive} logic frame the syntax and
the set of axioms are recursively coded. A recursive logic frame
is called recursively (countably) compact, if every recursive
(respectively, countable) finitely consistent theory has a model.
We show that  for logic frames built from the cardinality
quantifiers "there exists at least \(\lambda\)" recursive
compactness always implies countable compactness. On the other
hand we show that a recursively compact extension need not be
countably compact.
\end{abstract}



\section{Introduction}

For the definition of an abstract logic and  a generalized
quantifier the reader is refereed to \cite{MR87g:03033},
\cite{MR39:5329},  and \cite{MR39:5330}.
%
Undoubtedly the most important among abstract logics are the ones
that have a complete axiomatization of validity. In many cases,
most notably when we combine even the simplest generalized
quantifiers, completeness of an axiomatization cannot be proved in
ZFC alone but depends of principles like CH or $\diamondsuit$. Our
approach is look for ZFC-provable relationships between
completeness, recursive compactness and countable compactness of a
logic, that would reveal important features of the logic even if
we cannot settle any one of these properties per se. For example,
the countable compactness of the logic $\fol(Q_1,Q_2,Q_3,\ldots)$
cannot be decided in ZFC, but we prove in ZFC that if this logic
is recursively compact, it is countably compact. We show by
example that recursive compactness does not in general imply
countable compactness.

Examples of logics that have a complete axiomatization at least
under additional set theoretic assumptions are:

\begin{itemize}

\item The infinitary language \(L_{\omega_1\omega}\mbox{
\cite{MR31:1178}}\).

\item Logic with the generalized quantifier
$$Q_\alpha x\phi(x,\vec{y})
  \iff |\{x:\phi(x,\vec{y})|\ge\aleph_\alpha
  \mbox{ \cite{MR29:3364}}.$$

\item Logic with the {\it Magidor-Malitz quantifier}
$$Q^{\mbox{{\tiny MM}}}_\alpha
  xy\phi(x,y,\vec{z})\iff
  \exists X(|X|\ge\aleph_\alpha\wedge
  \forall x,y\in X\phi(x,y,\vec{z}))
  \mbox{ \cite{MR56:11746}}.$$

\item Logic with the {\it cofinality quantifier}
$$Q^{\mbox{{\tiny cof}}}_{\al_0}xy\phi(x,y,\vec{z})\iff
\{\la x,y\ra:\phi(x,y,\vec{z})\}\mbox{ has
cofinality \(\aleph_0\) \cite{MR51:12510}}.$$
\end{itemize}

\section{Logic Frames}

Our concept of a logic frame captures the combination of syntax,
semantics and proof theory of an extension of first order logic.
This is a very general concept and is not defined here with
mathematical exactness, as we do not prove any general results
about logic frames. All our results are about concrete examples.

\begin{definition}
\begin{enumerate}
\item
A {\it logic frame} is a triple
\[L^*=\la\L,\models_\L,\A\ra,\]
where

\begin{enumerate}

\item \(\la\L,\models_\L\ra\) is a logic in the sense of
Definition~1.1.1 in \cite{MR87g:03033}.

\item \(\A\) is a class of
\(L^*\)-axioms
and \(L^*\)-inference rules.

\end{enumerate}
\end{enumerate}
\end{definition}

We write
\(\vdash_\A\phi\) if \(\phi\)
is derivable using the axioms and rules in \(\A\).

\begin{example}\label{infinitary1}
Let
\[L_{\kappa\lambda}=
\la\L_{\kappa\lambda}, \models_{L_{\kappa\lambda}},
\A_{\kappa\lambda}\ra,\] where \(\A_{\kappa\lambda}\) has as
axioms the obvious axioms and Chang's Distributive Laws, and as
rules Modus Ponens, Conjunction Rule, Generalization Rule and the
Rule of Dependent Choices from \cite{MR31:1178}. This an old
example of a logic frame introduced by Tarski in the late 50's and
studied intensively, e.g. by Karp \cite{MR31:1178}.

\end{example}

\begin{example}\label{gq1}
Let
\[L(Q_\aa)=
\la\L(Q_\aa),
\models_{L(Q_\aa)},
\A_{Q_\aa}\ra,\]
where \(\A_{Q_\aa}\) has as axioms
the basic axioms of first order logic
and
\begin{eqnarray*}
&&\neg Q_\aa(x=y\vee x=z)\\
&&\forall x(\phi\ri\psi)\ri(Q_\aa x\phi\ri Q_\aa x\psi)\\
&&Q_\aa x\phi(x)\leftrightarrow Q_\aa y\psi(y),\mbox{ where
  \(\phi(x,...)\) is a formula of }\\
&&\mbox{\hspace{2cm}\(L(Q_\aa)\) in which
  \(y\) does not occur}\\
&&Q_\aa y\exists x\phi\ri\exists x Q_\aa y\phi\vee
     Q_\aa x\exists y\phi,
\end{eqnarray*}
 and Modus Ponens as the only rule.
The logic $L(Q_\aa)$  was introduced by Mostowski \cite{MR19:724f}
and the above frame by   Keisler \cite{MR41:8217}.
\end{example}

\begin{example}\label{MM1}
{Magidor-Malitz quantifier}
logic frame is
\[L(\qmm_\alpha)=
\la\L(\qmm_\alpha), \models, \A^{\mbox{\tiny MM}}_\alpha\ra,\]
where
$$Q^{\mbox{{\tiny MM}}}_\alpha
xy\phi(x,y,\vec{z})\iff \exists X(|X|\ge\aleph_\alpha\wedge
X\times X\subseteq \{\la x,y\ra:\phi(x,y\vec{z})\}) \mbox{
\cite{MR56:11746}}$$ and \(\A^{\mbox{\tiny MM}}_\alpha\) is the
set of axioms and rules introduced by Magidor and Malitz in
\cite{MR56:11746}.
\end{example}

\begin{definition}
\begin{enumerate}

\item A logic frame \(L^*=(\L,\models_\L,\A)\) is
{\em recursive} if

\begin{enumerate}

\item There is an effective algorithm which gives for each finite
vocabulary \(\tau\) the set \(\L[\tau]\) and for each \(\phi\in
\L[\tau]\) a second order formula\footnote{Second order logic
represents a strong logic with an effectively defined syntax. It
is not essential, which logic is used here as long as it is
powerful enough.} which defines the semantics of \(\phi\).

\item There is an effective algorithm
which gives the axioms and rules of \(\A\).

\end{enumerate}

\item A logic frame \(L^*=\la\L,\models_\L,\A\ra\) is
a {\em \(\la\kappa,\lambda\ra\)-logic frame}, if
each sentence contains less than \(\lambda\)
predicate, function and constant symbols,
and \(|\L[\tau]|\le\kappa\) whenever
the vocabulary \(\tau\) has less that \(\lambda\)
symbols altogether.

\end{enumerate}

\end{definition}

\begin{example}
The logic frame \(L_{\kappa\lambda}\) is a
\(\la\kappa^\kappa,\kappa\ra\)-logic frame. It is effective, if
\(\pk=\pl=\omega\). \label{gq2} The logic frame \(L(Q_\aa)\) is an
effective \(\la\omega,\omega\ra\)-logic frame. The logic frame
\(L(Q^{\mbox{\tiny MM}}_\alpha)\) is an effective
\(\la\omega,\omega\ra\)-logic frame.
\end{example}

\begin{definition}
A logic frame \(L^*=\la\L,\models_\L,\A\ra\)
is:
\begin{enumerate}

\item  {\em complete}
if every  \(\A\)-consistent
\(\L\)-sentence has a model.

\item {\em recursively compact} if every \(L^*\)-theory which is
recursive in the set of axioms and rules, and which has the
property that every finite subset of it has a model, has itself a
model.

\item {\em \((\kappa,\lambda)\)-compact}
if every  \(L^*\)-theory of cardinality
\(\le\kappa\), every subset of cardinality
\(<\lambda\) of which is \(\A\)-consistent,
has a model.

\item {\em countably compact}, if it is
\((\omega,\omega)\)-compact.
\end{enumerate}

\end{definition}

For recursive logic frames recursive compactness has a simpler
definition: Every recursive theory, every finite subset of which
has a model, has itself a model.

The logic frame \(L_{\kappa\lambda}\)
is  complete if
\begin{eqnarray*}
&1.&\kappa=\mu^+\mbox{ and } \mu^{<\lambda}=\mu\mbox{, or}\\
&2.&\kappa\mbox{ is strongly inaccessible, or}\\
&3.&\kappa\mbox{ is weakly inaccessible, }
  \lambda\mbox{ is regular and}\\
&&\hspace{1cm}(\forall\alpha<\kappa)
  (\forall\beta<\lambda)
  (\alpha^{\beta}<\kappa)\cite{MR31:1178}.\\
\end{eqnarray*}
\(L_{\kappa\lambda}\) does not satisfy the completeness theorem if
\(\kappa=\lambda\) is a successor cardinal (D.Scott, see
\cite{MR31:1178}). \(L_{\kappa\lambda}\) is not
\((\kappa,\kappa)\)-compact unless \(\kappa\) is weakly compact,
and then also \(L_{\kappa\kappa}\) is \((\kappa,\kappa)\)-compact.
\(L_{\kappa\lambda}\) is not strongly compact unless \(\kappa\) is
and then also \(L_{\kappa\kappa}\) is. If
\(\al_\aa^{<\al_\aa}=\al_\aa\), then by Chang's Two-Cardinal
Theorem, \(L(Q_{\aa+1})\) is complete and countably compact. If
V=L, then \(L(Q_{\aa})\) is complete and countably compact for all
\(\aa\). The logic frame \(L(Q^{\mbox{\tiny MM}}_\alpha)\) is
complete, if we assume \(\Diamond\), \(\Diamond_\aa\) and
\(\Diamond_{\aa+1}\), but there is a forcing extension in which
\(L(Q^{\mbox{\tiny MM}}_1)\) is not countably compact
\cite{AbSh:403}.

Completeness, which can always be achieved by adding new axioms,
does not imply recursive or countable compactness. However, if the
axioms and rules have a "finite character", as is the case in
first order logic, then the implication is true. Likewise,
recursive compactness does not, a priori, imply countable
compactness (see below for an example), although usually
counter-examples to compactness in extensions of first order logic
are very simple theories. This motivates the following definition:

\begin{definition}
A logic frame \(L^*=\la\L,\models_\L,\A\ra\) has

\begin{enumerate}

\item {\em finite recursive character} if for every possible
universe\footnote{I.e. inner model of forcing extension.} \(V'\)
\[V'\models(\mbox{\(L^*\) is complete
\(\Rightarrow\) \(L^*\) is recursively compact}).\]

\item {\em finite
\((\kappa,\lambda)\)-character} if for
every possible universe \(V'\)
\[V'\models(\mbox{\(L^*\) is complete
\(\Rightarrow\) \(L^*\) is
\((\kappa,\lambda)\)-compact}).\]

\item {\em recursive
\((\kappa,\lambda)\)-character} if for
every possible universe \(V'\)
\[V'\models(\mbox{\(L^*\) is recursively compact
\(\Rightarrow\) \(L^*\) is
\((\kappa,\lambda)\)-compact}).\]

\end{enumerate}

Mere "character" means \((\omega,\omega)\)-character. "Strong
character", means \((\kappa,\omega)\)-character for all
\(\kappa\).

\end{definition}

The definition of logic frames leaves many details vague, e.g. the
exact form of axioms and rules. Also the conditions of a recursive
logic frame would have to be formulated more exactly for any
general results. Going into such details would take us too much
astray from the main purpose of this paper.

\begin{example}
The logic frame \(L_{\kappa\lambda}\) is not of finite
\((\kappa,\kappa)\)-character, unless \(\kappa=\omega\), since it
is in some possible universes complete, but not
\((\kappa,\kappa)\)-compact. The logic frames \(L(Q_\aa)\) and
\(L(Q^{\mbox{\tiny MM}}_\aa)\) are in some possible universes
complete, but in some not countably compact. We discuss below the
problem whether they have recursive or finite character.

\end{example}


\section{A logic with  recursive character}


Let us consider the logic
\[\L=\fol(Q_{\aa_n})_{n<\oo},\]
where \(0<\aa_0<\aa_1<\ldots\) are arbitrary ordinals. We cannot
say in general whether \(\L\) is countably compact or not. If CH
holds, then $\L$ is countable compact \cite{Sh:8}, but it is
consistent that $\fol(Q_{n+1})_{n<\oo}$ is not countably compact
\cite{Shelah604}. There is a natural axiom system $\A$ for \(\L\)
based on so called identities. Using the methods of \cite{Sh:8} it
follows that if \(\al_0\) is small for each \(\al_{\aa_n}\)
(\(\mu\) is small for \(\ll\) if for every \(\ll_i\), \(i<\mu\),
we have \(\prod_{i<\mu}\ll_i<\ll\)), then this axiom system is
complete and \(\L\) is countably compact. In this section we show
that if \(\L\) is recursively compact, then \(\L\) is countably
compact. Thus \(\L\) gives rise to an example of a logic frame
with recursive character.

The model theory of \(\L\) is closely tied with the model theory
of \((\al_{\aa_n})_{n<\oo}\)-like models. This follows from usual
reduction techniques (see \cite[p. 45]{MR87g:03033}). We define
now a generalized concept of identity needed for the formulation
of the axioms of the logic \(\L\). The concept of identity was
introduced in \cite{Sh:8} (see also \cite[p. 188]{MR87g:03033} for
a survey). A generalized identity is a sequence of finite
equivalence relations on finite sets such that equivalent sets
have the same size. In addition, the generalized identities are
attached with a finite sequence of functions.

\begin{definition}
\begin{enumerate}
\item A {\em} generalized identity is
a triple
\[\I=\II,\]
where
\begin{enumerate}
\item Each \(D_n\) is a finite set of ordinals,

\item \(n< n'\) implies \(D_n\subseteq D_{n'}=\emptyset\),

\item Each \(I_{n,a}\), $a\in D_n$,  is an equivalence relation on
\(D_n^{<\oo}\),

\item \(xI_{n,a}y\) implies \(|x|=|y|\) for all \(x\) and \(y\),

\item \(h_n:D_n^{<\oo}\rightarrow D_n\) and $n<m$ implies
$h_n\subseteq h_m$.
\end{enumerate}

\item Suppose \[\I=\II\] and \[\I'=\IIP\] are generalized
identities. We say that \(\I\) is a {\em subidentity}  of \(\I'\),
if there are a one-to-one \(\ps:m\ri m'\) and an order-preserving
\(\pi:D_i\ri D'_{\ps(i)}\) such that \[xI_{n,a} y\iff
\pi[x]I'_{\ps(n),\pi(a)} \pi[y]\] and \(\pi(h_i(x))\ge
h'_i(\pi[x])\). If such bijections \(\pi\) and \(\ps\) exist, then
the generalized identities are called {\em equivalent}.

\item Let \(\oF\) be the class of all \(\F=\la
(F_{l,\aa})_{l<\oo,\aa<\oo_{\aa_l}},(h_l)_{l<\oo}\ra\), where
\[F_{l,\aa}:[\oo_{\aa_l}]^{<\oo}\ri\kk_{l,\aa}\] such that
\(\al_0\le \kk_{l,\aa}<\al_{\aa_l}\) and
\[h_l:[\oo_{\aa_l}]^{<\oo}\ri\oo_{\aa_l}.\]

\item If \(X\se\bigcup_{l<\oo}\oo_{\aa_l}\) is finite, then \(\F\)
{\em induces} the generalized identity
\[\I=\II,\] where
\begin{enumerate}
\item \(\{l:\al_{\aa_l}\in X\}=\{l_0,...,l_{m-1}\}\),

\item \(D_n=\oo_{\aa_{l_n}}\cap X\),

\item \(xI_{i,\aa} y\) if and only if
\(F_{l_i,\aa}(x)=F_{l_i,\aa}(y)\).
\end{enumerate}

\item \(\I(\F)\) denotes the set of identities which are
subidentities of those induced by \(\F\).
$\Ifin(\F)=\{\I\in\I(\F):\forall n<m(D_n\subseteq\oo)\}$

\item \(\I(\vecz{\aa})=\bigcap\{\I(\F):\F\in\oF\}\), where
$\vecz{\aa}=\la\aa_n:n<\oo\ra$.

\item A {\em fundamental function} for $\vecz{\aa}$ is an $\F$
such that $\I(\F)=\I(\vecz{\aa}).$
\end{enumerate}
\end{definition}

A priori, the elements of the domains $D_n$ of a generalized
identity can be any ordinals, but up to equivalence, they can
always be taken to be natural numbers.  Thus, if we assume a
canonical coding of such generalized identities by natural
numbers, it makes sense to ask whether a certain set of identities
is recursive or not. Also, $\Ifin(\F)=\I(\F)$ up to equivalence.

\begin{definition}
Suppose \(\aa_n, n<\oo\) are ordinals in increasing order. An
\((\lan)_{n<\oo}\)-like model is a model \(\ma\) in a language
with distinguished predicates \(P_n, n<\oo,\) and a binary
predicate \(<\) such that for all \(n\) \(\la P_n^\ma,<^\ma\ra\)
is \(\lan\)-like, i.e. a linear order of cardinality \(\lan\) and
every initial segment is of cardinality \(<\lan\).
\end{definition}

The definition of a generalized identity looks complicated but its
meaning becomes completely transparent when one realizes that it
is exactly what one needs to construct by means of Skolem
functions and the Compactness Theorem an \((\lan)_{n<\oo}\)-like
model for a first order theory.

\begin{proposition}\label{compakt} If there is a fundamental function for
$\vecz{\aa}$, then $\fol(Q_{\aa_n})_{n<\oo}$ is
$({\ll},\oo)$-compact for any $\ll<\al_{\aa_0}$.
\end{proposition}

\proof This is like \cite{Sh:8}. Suppose \(\F=\la
(F_{l,\aa})_{l<\oo,\aa<\oo_{\aa_l}},(h_l)_{l<\oo}\ra\) is a
fundamental function for $\vecz{\aa}$. Suppose $T$ is a finitely
consistent $\fol(Q_{\aa_n})_{n<\oo}$ theory. W.l.o.g. $T$ has
built-in Skolem functions, the language of $T$ includes unary
predicates $P_l$, $l<\oo$, and a binary predicate $<$, and it
suffices to construct a $(\la\lal)_{l<\oo}$-like model for $T$.
Let $T^*$ consist of $T$ plus the axioms

\begin{enumerate}

\item $c^l_\xi<c^l_\eta$ for $\xi<\eta<\lal$ and $l<\oo$,

\item $P_l(c^l_\xi)$ for $\xi<\lal$ and $l<\oo$,

\item
$t(c^l_{\xi_1},\ldots,c^l_{\xi_n})<c^l_{h_l(\xi_1,\ldots,\xi_n)}$
 for all terms $t$ and $(\xi_1,\ldots,\xi_n)\in (\lal)^n$,

\item $\begin{array}[t]{l}
t(c^l_{\xi_1},\ldots,c^l_{\xi_n})=t(c^l_{\eta_1},\ldots,c^l_{\eta_n})\\
\vee(c^l_\alpha\le t(c^l_{\xi_1},\ldots,c^l_{\xi_n})\wedge
c^l_\alpha\le t(c^l_{\eta_1},\ldots,c^l_{\eta_n}))\end{array}$

 for all terms $t$ and $(\xi_1,\ldots,\xi_n)\in (\lal)^n$ such
 that $F_{l,\aa}(\xi_1,\ldots,\xi_n)=F_{l,\aa}(\eta_1,\ldots,\eta_n)$
\end{enumerate}

\noindent It is easy to see that the Skolem closure of the
constants is  $(\la\lal)_{l<\oo}$-like in every model of $T^*$.
Thus it suffices to show that $T^*$ is finitely consistent. We
refer to the proof of Theorem 3.2.1 in \cite[Chapter
V]{MR87g:03033} for details. \qed

A consequence of the above proof is:

\begin{corollary} There is a set $Ax(\vecz{\aa})$ of valid sentences of
$\fol(Q_{\aa_n})_{n<\oo}$ such that if there is a fundamental
function for $\vecz{\aa}$, then a sentence $\phi$ of
$\fol(Q_{\aa_n})_{n<\oo}$ is valid if and only if it follows from
$Ax(\vecz{\aa})$ and the axiom schemas of first order logic using
rules of proof of first order logic. If {\I(\vecz{\aa})} is r.e.,
then so is $Ax(\vecz{\aa})$.
\end{corollary}

The axiomatization $Ax(\vecz{\aa})$ may not be complete, but the
point is, that {if} there is a fundamental function for
$\vecz{\aa}$, {then} it is complete. So we know the
axiomatization, but we do not always have a fundamental function.
Likewise, we do not know in general whether $\I(\vecz{\aa})$ is
r.e. but if it is, then $Ax(\vecz{\aa})$ gives a recursive
complete axiomatization of $\fol(Q_{\aa_n})_{n<\oo}$. Thus the
triple
$$\fol(Q_{\aa_n})_{n<\oo}=\la {\cal
L}(Q_{\aa_n})_{n<\oo},\models,Ax(\vecz{\aa})\ra$$ forms a logic
frame, which is a recursive logic frame if $\I(\vecz{\aa})$ is
r.e. and complete if there is a fundamental function for
$\vecz{\aa}$.

\begin{lemma}
Suppose \(\vecz{\aa}=\la\aa_n: n<\oo\ra\) is a sequence of
ordinals in increasing order. Suppose \(\I\) is a generalized
identity. There is a sentence \(\phi_\I\) in
\(L(Q_{\aa_n})_{n<\oo}\) such that the following conditions are
equivalent:
\begin{enumerate}
\item \(\phi_\I\) has a model.

\item \(\I\notin\I(\vecz{\aa})\)
\end{enumerate}
\end{lemma}

\proof Suppose $\I=\II$. The vocabulary of \(\phi_{n,\I}\) has a
unary predicate \(P_l\), a binary predicate \(<\) and an \(i\)-ary
function symbol \(F^i_{l,a}\) for each \(l<m\), $a\in D_l$ and
\(i<n=|\bigcup_{j<m}D_j|\). Let \(\phi^-_\I\) be the conjunction
of
\begin{enumerate}
\item \(\la P_l,<_l\ra\) is a \(\lal\)-like linear order for
$l<m$,

\item \(F^i_{l}:P_l\times (P_l)^i \to P_k\),

\item $F^i_{l}(a,x)<g_{l}(a)$,

\item \(h^i_{l}: (P_l)^i \to P_l\).

\end{enumerate}
Any model $M$ of $\phi^-_\I$ givers rise to \(\F_M=\la
(F_{l,\aa})_{l<\oo,\aa<\oo_{\aa_l}},(h_l)_{l<\oo}\ra\), where
$$F_{l,a}(x)=(F^i_{l,a})^M(x), \mbox{ if $x\in ((P_l)^M)^i$}$$
 and
$$h_{l}(x)=(h^i_{l})^M(x), \mbox{ if $x\in ((P_l)^M)^i$}.$$
The sentence \(\phi_\I\) is the conjunction of \(\phi^-_\I\) and a
first order sentence stating that in models $M$ of \(\phi_\I\) the
identity  $\I$ is not a subidentity of the identity induced by
$\F_M$. \qed

\begin{lemma}\label{ra}If the $\fol(Q_{\aa_n})_{n<\oo}$-theory
$\{\sigma_\I:\I\notin\I(\vecz{\aa})\}$ has a model, there is a
fundamental function for $\vecz{\aa}$.
\end{lemma}

\proof Suppose $M$ is such a model. Let $\F_M$ be as in the
previous proof. It is clear that $\F_M$ is a fundamental function.
\qed

\begin{lemma}\label{ha}Every finite subset of  the
$\fol(Q_{\aa_n})_{n<\oo}$-theory
$\{\sigma_\I:\I\notin\I(\vecz{\aa})\}$ has a model.
\end{lemma}

\proof Suppose $\I\notin\I(\vecz{\aa})$. It suffices to show that
$\sigma_\I$ has a model. Let $\F=\la
(F_{l,\aa})_{l<\oo,\aa<\oo_{\aa_l}},(h_l)_{l<\oo}\ra$ be such that
$\I$ is not a subidentity of one induced by $\F$. It is easy to
build a $(\la\lal)_{l<\oo}$-like model $M$ from $\F$ so that
$M\models\sigma_\I$. \qed

\begin{theorem} The logic frame  $\fol(Q_{\aa_n})_{n<\oo}$ has recursive character
for all $\vecz{\aa}$.
\end{theorem}

\proof If there is a fundamental function for $\vecz{\aa}$, then
$\fol(Q_{\aa_n})_{n<\oo}$ is countably compact by
Proposition~\ref{compakt}. Thus we may assume that there is no
fundamental function for $\vecz{\aa}$. Let $\J=\{\I_n:n<\oo\}$ be
a canonical enumeration of all $\I\notin\I(\vecz{\aa})$. Let $T$
be the $\fol(Q_{\aa_n})_{n<\oo}$-theory consisting of
\begin{enumerate}

\item "$c$ has at least $n$ predecessors in $<$", for $n<\oo$,

\item "If $c$ has at least $n$ predecessors in $<$, then
$\sigma_{\I_n}$".
\end{enumerate}

\noindent The theory $T$ is recursive in $\J$, and  finitely
consistent by Lemma~\ref{ha}. On the other hand, any model of $T$
would give rise to a fundamental function for $\vecz{\aa}$ by
Lemma~\ref{ra}. \qed


\section{A logic which does not have recursive character}


We show that there is a logic frame $L^*$ which is recursively
compact but not countably compact.
 We make use of  the  quantifier \(\qst\) from
\cite{MR51:12510}. To recall the  definition of \(\qst\) we adopt
the following notation:

\begin{definition}\label{filtration}
Let \(\ma=(A,R)\) be an arbitrary \(\aleph_1\)-like linearly
ordered structure. We use $H(\ma)$ to denote the set of all
initial segments of $\ma$. A {\em filtration} of \(\ma\) is subset
\(X\) of \(H(\ma)\) such that $A=\bigcup_{I\in X}I$ and $X$ is
closed under unions of increasing sequences. Let $D(\ma)$ be the
filter on $H(\ma)$ generated by all filtrations of $\ma$.
\end{definition}

The quantifier \(\qst\) is defined as follows:

\begin{definition}\label{stationary logic as lindstrom}
The generalized quantifier $\qst$ is defined by \[\ma\models \qst
xy\phi(x,y,\vec{a})\] if and only if
\[\ma=(A, R_\phi),\mbox{ where }R_\phi=\{(b,c):\ma\models\phi(b,c,\vec{a})\}\]
 is an $\aleph_1$-like linear
ordered structure  such that \[\{I\in H(\ma):I \mbox{ does not
have a sup in }R_\phi\}\notin D(\ma).\]
\end{definition}

It follows from \cite{Sh:43} and \cite{MR82f:03031a} that
$\fol(\qst)$ equipped with some natural axioms and rules is a
complete countably compact logic frame.

\begin{definition} If $S\subseteq\omega_1$, then the generalized
quantifier $\qsts{S}$
 is defined by
\[\ma\models \qsts{S} xy\phi(x,y,\vec{a})\] if and only if
\(R_\phi\) is an $\aleph_1$-like linear order of \(A\) with a
filtration \(\{ I_\alpha:\alpha<\omega_1\}\) such that
\[\forall\alpha<\omega_1((I_\alpha \mbox{
has a sup in }R_\phi)\iff \alpha\in S).\] The syntax of the logic
\[\st\] is defined as follows: \(\st\) extends first order logic
by the quantifiers $Q_1$, \(\qst\) and the infinite number of new
formal quantifiers \(\qsts{X_n}\) (we leave $X_n$
unspecified).\end{definition}

If we fix a sequence $\ss$ and let $\qsts{X_n}$ be interpreted as
$\qsts{S_n}$, we get a definition of semantics of $\st$. We call
this semantics the {\em $\ss$-interpretation} of $\st$.


\begin{definition} We call a finite sequence \(\sigma=\ssn\)
(or an infinite sequence \(\ss\)) of subsets of \(\omega_1\) {\em
stationary independent}, if all finite Boolean combinations of the
sets \(S_i\) are stationary.
\end{definition}

We will show now that the set of valid sentences of $\st$ is
independent of the sequence $\sss$, as long as this sequence is
stationary independent.

Rather than giving an explicit axiom system for \(\st\) we
manipulate models of set theory in order to get the same  results.
Throughout, we use \(\zfc\) to denote a finite fragment of
\(\mbox{ZFC}\) sufficient for the arguments involved. This is only
to guarantee the existence of models.  We make extensive use of
the following result of Hutchinson:

\begin{lem}\cite{MR53:12944}\label{hutchinson}
Suppose \(\mm\) is a countable model of \(\zfc\) and \(S\in M\)
such that
\[\mm\models "S\subseteq\omega_1\mbox{ is stationary}".\]
Then there is a countable  \(\mn\) such that \(\mm\prec\mn\),
\(\mn\) has a least new ordinal \(\delta\), \(\mn\models \delta\in
S\), and every \(a\in M\) such that \(\mm\models
"a\epsilon\omega_1"\) remains fixed i.e. \[\forall b\in
N([\mn\models b\epsilon a] \rightarrow b\in M).\]
\end{lem}
 As
pointed out in \cite[Theorem 3.5]{MR55:10268},
Lemma~\ref{hutchinson} can be iterated \(\omega_1\) times to get
an  elementary extension \(\mn\) of \(\mm\) such that
\(\la\omega_1^\mn,\epsilon^\mn\ra\) is an \(\aleph_1\)-like linear
order with a filtration \(\la C_\alpha:\alpha<\omega_1\ra\) such
that \(C_\alpha\in N\) and
\[\{\alpha:\mn\models"C_\alpha\mbox{ has a sup in }S"\}\]
is a club subset of \(\omega_1\). In the proof of \cite[Theorem
3.2.2]{MR87g:03033} this idea is elaborated by splitting
\(\omega_1\) first into \(\aleph_1\) disjoint stationary sets and
then extending all possible stationary sets, one at a time. This
yields an elementary extension \(\mn\) of \(\mm\) such that
\begin{itemize}
\item[(G)] \(\la\omega_1^\mn,\epsilon^\mn\ra\) is an
\(\aleph_1\)-like linear order with a filtration \(\la
C_\alpha:\alpha<\omega_1\ra\) such that each \(C_\alpha\) is in
\(N\) and if \(\mn\models"S\subseteq\omega_1\mbox{ is
 stationary}"\), then
\[\{\alpha:\mn\models"C_\alpha\mbox{ has a sup in }S"\}\]
is a stationary subset of \(\omega_1\).
\end{itemize}
Let us call a model \(\mm\) of set theory {\it good}, if it
satisfies (G). We have sketched a proof of:
\begin{lem}
Every countable  model \(\mm\) of \(ZFC\) has a good
 elementary extension of cardinality
\(\aleph_1\).
\end{lem}
With a minor modification we get:
\begin{lem}\label{lhsu}
Suppose \(\mm\) is a countable model of \(\zfc\) and \(\sigma\in
M\) such that
\[\mm\models "\sigma=\ss\mbox{ is stationary independent}."\]
Suppose \(\ttn\) is an arbitrary stationary independent sequence.
Then there is a good elementary extension \(\mn\) of \(\mm\) such
that \(\la\omega_1^\mn,\epsilon^\mn\ra\) is an \(\aleph_1\)-like
linear order with a filtration \(\la C_\alpha:\alpha<\omega_1\ra\)
and a club \(D\) such that for all \(\alpha<\omega_1\)
\(C_\alpha\in N\) and for all \(i=0,...,n\) and all \(\alpha\in
D\) we have
\[\mn\models"C_\alpha\mbox{ has a sup in }S_i"\iff\alpha\in T_i.\]
\end{lem}

\proof One can imitate the proof in \cite[Theorem 3.5]{MR55:10268}
and the proof of \cite[Theorem 3.2.2]{MR87g:03033}.
W.l.o.g., the sets in \(\ssn\) partition
\(\omega_1\) in \(\mm\)
and the same holds for \(\ttn\).
In the iteration of Lemma~\ref{hutchinson} we extend at
stage \(\xi\) the set \(S_i\) if \(\xi\in T_i\). \qed

%
%

\begin{pro}\label{whatever} Suppose $\phi$ is a sentence of $\st$. Then the
following conditions are equivalent:
\begin{enumerate}
\item $\phi$ has a model in some $\ss$-interpretation, where $\ss$
is stationary independent.

\item $\phi$ has a model in all $\ss$-interpretations, for $\ss$
stationary independent.
\end{enumerate}
\end{pro}

\proof  Suppose \(\phi\in\st\). Choose \(n\in \omega\) such that
\(\phi\) contains no quantifiers \(\qsts{X_i}\) for \(i > n\). Let
\(\Phi(\ma,\sigma,\phi)\) be a formula of set theory expressing
the conjunction of "\(\ma\models\phi\) in the
$\sigma$-interpretation" and "\(\sigma=\ssn\) is stationary
independent" in such a way that if \(\mm\) is a good model of
\(\zfc\) containing \(\ma\) and \(\sigma\), then the following
conditions are equivalent:
\begin{itemize}
\item \(\la\ma,\sigma\ra\models\phi\). \item
\(\mm\models\Phi(\ma,\sigma,\phi)\)
\end{itemize}
The claim follows now from Lemma \ref{lhsu}. \qed

\begin{dfn} Let $\Val(\st)$ be the set of sentences of
$\st$ which are valid under  $\ss$-interpretation for some
(equivalently, all) stationary independent $\ss$.
\end{dfn}

\begin{pro}
The set $\Val(\st)$ is recursively enumerable, provably in ZFC.
\end{pro}

\proof Suppose \(\phi\in\st\). Let $\Phi$ be as in the proof of
Proposition~\ref{whatever}. By Lemma~\ref{lhsu} we have the
equivalence of
\begin{itemize}
\item \(\neg\phi\notin \Val(\st)\). \item
\(\zfc\cup\{\exists\ma\exists\sigma\Phi(\ma,\sigma,\phi)\}\)
      is consistent in $\fol$.
\end{itemize}
Since the latter is a \(\Pi^0_1\)-property of \(\phi\), we have
proved the claim. \qed

\begin{cor} There is a recursive set $\Ax(\st)$ of sentences of $\st$ such that a
sentence of $\st$ is in $\Val(\st)$ if and only if it follows from
$Ax(\st)$ and the axiom schemas of first order logic using rules
of proof of first order logic.
\end{cor}

\begin{dfn} Suppose $\ss$ is stationary independent. We define a new  recursive logic frame
$$\st(\sss)=\la\st,\models,{\cal A}\ra,$$ where ${\cal A}$ consists of
$Ax(\st)$ and the axioms and rules of first order logic.
\end{dfn}

\begin{cor}[Completeness Theorem for $\st(\sss)$] \label{taydelliosyyslause}
The logic frame $\st(\sss)$ is complete.
\end{cor}

The axioms of $\st(\sss)$ state the stationary independence of
$\ss$. Thus the completeness of $\st(\sss)$ is vacuous for
stationary nonindependent $\ss$. The same method gives the
following partial countable compactness result: Suppose $\sss$ is
stationary independent. Any finitely consistent countable theory
in $\st$, which contains an occurrence of $\qsts{X_n}$ for only
finitely many $n$, has a model.

Note that the syntax and the axioms of $\st(\sss)$ are independent
of $\ss$. We conjecture that there is a natural complete axiom
system for all $\st(\sss)$ based on
\begin{itemize}

\item The usual axioms and rules of $\L(Q_1)$ as in
\cite{MR41:8217}.

\item Natural axioms (like Fodor's Lemma) for $\qst$ as in
\cite{MR82f:03031a}.

\item Axiom schemas stating the stationary independence of $\ss$.
\end{itemize}

Let $\sti(\sss)$ be the extension of $\st$ obtained by allowing
countable conjunctions and disjunctions.

\begin{pro}\label{complete2}
Suppose \(\phi\in\sti\). The predicate "\(\phi\mbox{ has a
model}\)" is
 a \(\Sigma_1^{ZFC}\)-definable property of
\(\phi\).
\end{pro}

\proof It suffices to notice that if $\mm$ is $\omega$-standard in
Lemma~\ref{hutchinson}, then so is $\mn$. \qed

 By making different choices for the stationary independent $\ss$, we can get
logics with different properties. We illustrate this now by making
a choice of $\ss$ which will render $\st(\sss)$ recursively
compact but not countably compact.

Let us fix a countable vocabulary \(\tau\) which contains
infinitely many symbols of all arities. Let \(T_n\), \(n<\omega\),
list all $Ax(\st)$-consistent recursive \(\st\)-theories in the
vocabulary \(\tau\). Let \(\tau^n\) be a new disjoint copy of
\(\tau\) for each \(n<\omega\). Let \(\tau^*\) consist of the
union of all the \(\tau_n\), the new binary predicate symbol
\(<^*\), and new unary predicate symbols \(P_n\) for \(n<\omega\).
If \(\phi\) is a formula and \(d\in 2\), let \((\phi)^d\) be
\(\phi\), if \(d=0\), and \(\neg\phi\), if \(d=1\). If
\(S\subseteq\omega_1\), then \((S)^{d}\) is defined similarly. For
any \(\eta:\omega\rightarrow 2\) let \(\psi_\eta\in\sti\) be the
conjunction of the following sentences of the vocabulary
\(\tau^*\):

\begin{enumerate}
\item[(a)] \(T_n\)
 translated into the vocabulary
\(\tau^n\).

\item[(b)] \(<^*\) is an $\aleph_1$-like linear order of the
universe.

 \item[(c)] \(\qsts{X_{n}} xy(x<^*y\wedge P_n(x)\wedge
P_n(y))\).

\item [(d)] \(\neg\exists x\bigwedge_{n}(P_n(x))^{\eta(n)}\).

\end{enumerate}

\begin{lem}
There is \(\eta:\omega\rightarrow 2\) such that $\psi_\eta$ has a
model.
\end{lem}

\proof Let $\Gamma$ consist of the sentences (a)-(c). By
Corollary~\ref{taydelliosyyslause}, \(\Gamma\) has a  model
\(\mm\) of cardinality \(\aleph_1\) in the $\ss$-interpretation
for some stationary independent $\ss$. Get a new
\(\eta:\omega\rightarrow 2\) by Cohen-forcing. Then in the
extension \(V[\eta]\)
\[
\bigcap_n (S_n)^{\eta(n)}=\emptyset.\]
 Thus \(V[\eta]\) satisfies the \(\Sigma_1\)-sentence
\begin{equation}\label{sigma1}
\exists\eta(\psi_\eta\mbox{ has a model}).
\end{equation}
By the Levy-Shoenfield Absoluteness Lemma and
Proposition~\ref{complete2} there is \(\eta\) in \(V\) such that
(\ref{sigma1}) holds in $V$. \qed

Now let $\sst$ be stationary independent such that $\psi_\eta$ has
a model $\mm^*$ in the $\sst$-interpretation.


\begin{theorem}
The recursive logic frame \(\st(\ssst)\) is recursively compact
but not countably compact.
\end{theorem}

\proof Suppose \(T\) is a consistent recursive theory in \(\st\).
W.l.o.g. \(T=T_m\) for some \(m<\omega\). Thus \(\mm^*\restriction
\tau^n\) gives immediately a model of \(T\). To prove that \(\st\)
is not countably compact, let \(T\) be a theory consisting of the
following sentences:
\begin{enumerate}
\item[(i)] \(<^*\) is an $\aleph_1$-like linear order.

\item[(ii)] \(\qsts{S^*_n} xy(x<^*y\wedge P_n(x)\wedge P_n(y))\)
for \(n<\omega\).

\item[(iii)] \(\qst xy(x<^*y\wedge P(x)\wedge  P(y))\).

\item[(iv)] \(\forall x(P(x)\rightarrow (P_n(x))^{\eta(n)})\) for
\(n<\omega\).
\end{enumerate}
Any finite subtheory of \(T\) contains only predicates
\(P_0,...,P_m\) for some \(m\), and has therefore a model: we let
\(P_i=S^*_i\) for \(i=0,...,m\) and
\[P=(P_0)^{\eta(0)}\cap...\cap (P_m)^{\eta(m)}.\] On the other
hand, suppose \(\la A,<^*,P,P_0,P_1,...\ra\models T\).
 By (ii) there are filtrations \(\la
D^n_\alpha:\alpha<\omega_1\ra\) of \(<^*\) and clubs \(E^n\) such
that for all \(n\) and for all \(\alpha\in E^n\)
\[\{\alpha<\omega_1:
D^n_\alpha\mbox{ has a sup in }\la A,<^*\ra\}=S^*_n.\]
By (iii) there is a
filtration
\(\la F_\alpha:\alpha<\omega_1\ra\) of \(<^*\) such that
\[B=\{\alpha<\omega_1:
F_\alpha\mbox{ has a sup in } P\}\] is stationary. Let
\(E^*\subseteq \bigcap_nE_n\) be a club such that
\(C_\alpha=D^n_\alpha=F_\alpha\) for \(\alpha\in E^*\) and
\(n<\omega\). Let \(\delta\in E^*\cap B\) and \(a=\sup F_\delta\).
Then \(a\in P\). Hence \(a\in\bigcap_n(P_n)^{\eta(n)}\) by (iv).
As \(a=\sup D_\delta^n\) for all \(n\), we have
\(a\in\bigcap_n(S^*_n)^{\eta(n)}\), contrary to the choice of
\(\eta\). We have proved that theory \(T\) has no models. \qed

Thus $\st$ does not have finite character. We end with an example
of a logic which, without being provably complete, has anyhow
finite character:

Recall that $\diamondsuit_S$ for $S\subseteq\oo_1$ is the
statement that there are sets $A_\aa\subseteq\aa$, $\aa\in S$,
such that for any $X\subseteq\omega_1$, the set $\{\aa\in S:X\cap
\aa=A_\aa\}$ is stationary.

\begin{dfn} Let $\L^{\diamondsuit}$ be the extension of $\fol$ by $Q_1$, $\qst$
and $\qsts{\sd}$, where

\[\sd=\left\{
\begin{array}{ll}
\emptyset,&\mbox{ if there is no bistationary $S$ with $\diamondsuit_S$}\\
\omega_1,&\mbox{ if there is a bistationary $S$ with $\diamondsuit_S$ but no maximal one}\\
S,&\mbox{ if $S$ is a maximal bistationary $S$ with $\diamondsuit_S$}\\
\end{array}
\right.\] We get a recursive logic frame
$L^{\diamondsuit}=\la\L^{\diamondsuit},\models,\A\ra$ by adapting
the set $\Ax(\st)$ to the case of just one bistationary set.
\end{dfn}

\begin{thm} $L^{\diamondsuit}$ has finite character.
\end{thm}

\proof Suppose there is no bistationary $S$ with $\diamondsuit_S$.
Then the consistent sentence ``$\mbox{$<$ is an $\aleph_1$-like
linear order}\wedge\qst
xy(x<y)\wedge\qsts{S_{\diamondsuit}}(x<y)$" has no model, so $\L$
is incomplete. Suppose there is a bistationary $S$ with
$\diamondsuit_S$ but no maximal one. Then the consistent sentence
``$\mbox{$<$ is an $\aleph_1$-like linear order}\wedge \qst
xy(x<y\wedge P(x))\wedge\qsts{S_{\diamondsuit}}(x<y\wedge \neg
P(x))$" has no model, so $\L$ is again incomplete. Finally,
suppose there is a maximal bistationary $S$ with $\diamondsuit_S$.
Now $\L$ is countably compact by the remark right after
Corollary~\ref{taydelliosyyslause}. \qed

Our results obviously do  not aim to be optimal. We merely want to
indicate that the concept of a logic frame offers a way out of the
plethora of independence results about generalized quantifiers.
The logic $\fol(Q_{n+1})_{n<\oo}$ is a good example. The results
about its countable compactness under CH and countable
incompactness in another model of set theory leave us perplexed
about the nature of the logic. Having recursive character reveals
something conclusive and positive, and raises the question, do
other problematic logics also have recursive character. Our logic
$\st$ is the other extreme: it is always completely axiomatizable,
but a judicious choice of $\ss$ renders it recursively compact
without being countably compact.
\bigskip

 \noindent{\bf Open
Question:} Does the Magidor-Malitz logic $L(\qmm_1)$ have
recursive character?




\def\cprime{$'$} \def\germ{\frak} \def\scr{\cal}
  \ifx\documentclass\undefinedcs\def\rm{\fam0\tenrm}\fi
  \def\defaultdefine#1#2{\expandafter\ifx\csname#1\endcsname\relax
  \expandafter\def\csname#1\endcsname{#2}\fi} \defaultdefine{Bbb}{\bf}
  \defaultdefine{frak}{\bf} \defaultdefine{mathbb}{\bf}
  \defaultdefine{mathcal}{\cal}
  \defaultdefine{beth}{BETH}\defaultdefine{cal}{\bf} \def\bbfI{{\Bbb I}}
  \def\mbox{\hbox} \def\text{\hbox} \def\om{\omega} \def\Cal#1{{\bf #1}}
  \def\pcf{pcf} \defaultdefine{cf}{cf} \defaultdefine{reals}{{\Bbb R}}
  \defaultdefine{real}{{\Bbb R}} \def\restriction{{|}} \def\club{CLUB}
  \def\w{\omega} \def\exist{\exists} \def\se{{\germ se}} \def\bb{{\bf b}}
  \def\equivalence{\equiv} \let\lt< \let\gt> \def\cite#1{[#1]}


{\tt shelah@math.huji.ac.il}

{\tt jouko.vaananen@helsinki.fi}
\end{document}